\documentclass[11pt,a4paper]{article}
\usepackage[utf8]{inputenc}
\usepackage{amsmath,amssymb,amsthm}
\usepackage{mathtools}
\usepackage{booktabs}
\usepackage{tabularx}
\usepackage{array}
\usepackage{microtype}
\usepackage{geometry}
\usepackage{algorithm}
\makeatletter
\providecommand{\algorithmname}{\fname@algorithm} 
\makeatother
\usepackage{algpseudocode}

\makeatletter
\newenvironment{breakablealgorithm}
  {%
    \begin{center}
      \refstepcounter{algorithm}%
      \hrule height.8pt depth0pt \kern2pt%
      \renewcommand{\caption}[2][\relax]{%
        {\raggedright\textbf{\algorithmname~\thealgorithm} ##2\par}%
        \ifx\relax##1\relax
          \addcontentsline{loa}{algorithm}{\protect\numberline{\thealgorithm}##2}%
        \else
          \addcontentsline{loa}{algorithm}{\protect\numberline{\thealgorithm}##1}%
        \fi
        \kern2pt\hrule\kern2pt
      }%
  }
  {%
      \kern2pt\hrule\relax%
    \end{center}
  }
\makeatother

\usepackage[authoryear,round]{natbib}
\setcitestyle{aysep={,},citesep={;},open={(},close={)}}

\usepackage[colorlinks=true,linkcolor=blue,citecolor=blue,urlcolor=blue]{hyperref}

\newtheorem{theorem}{Theorem}[section]
\newtheorem{lemma}[theorem]{Lemma}
\newtheorem{proposition}[theorem]{Proposition}

\theoremstyle{definition}
\newtheorem{definition}[theorem]{Definition}

\theoremstyle{remark}

\newcommand{\R}{\mathbb{R}}
\newcommand{\C}{\mathbb{C}}
\newcommand{\g}{\mathfrak{g}}
\newcommand{\n}{\mathfrak{n}}
\newcommand{\hh}{\mathfrak{h}}
\newcommand{\s}{\mathfrak{s}}
\newcommand{\rr}{\mathfrak{r}}
\newcommand{\ad}{\operatorname{ad}}
\newcommand{\Aut}{\operatorname{Aut}}
\newcommand{\Der}{\operatorname{Der}}

\newcommand{\tr}{\operatorname{tr}}
\newcommand{\diag}{\operatorname{diag}}

\title{\textbf{Classification of real Lie algebras of dimensions four}}
\author{Mehdi Nadjafikhah\thanks{School of Mathematics and Computer Sciences, Iran University of Science and Technology, Narmak, Tehran 16846-13114, Iran. email: m\_nadjafikhah@iust.ac.ir}}
\date{\today}

\begin{document}
\maketitle

\begin{abstract}
We present an exhaustive, coordinate-free, and self-contained classification of four -dimensional real Lie algebras up to isomorphism. In dimension three, classification schemes frequently exploit Bryant's classical exterior duality isomorphism $\Lambda^2\g \cong \g^* \otimes \Lambda^3\g$, which reduces the Jacobi identity to an elegant linear algebraic condition on a single $3\times 3$ matrix. We explicitly demonstrate that this accidental isomorphism breaks down in dimension four due to the dimensional mismatch $\dim \Lambda^2\g = 6 \neq 4 = \dim(\g^* \otimes \Lambda^4\g)$, transforming the Jacobi identity into a non-linear system of 16 quadratic constraints on 24 structure constants. 

To overcome this geometric obstruction without relying on brute-force polynomial elimination, we employ an intrinsic Lie-theoretic reduction: the Levi--Malcev decomposition isolates the non-solvable sector into direct products with three-dimensional simple ideals, while the solvable sector is classified by anchoring the structure to its nilradical $\n$, which is an intrinsic characteristic ideal. For codimension-one nilradicals ($\dim\n = 3$), the problem is reduced to classifying outer derivations of $\R^3$ and the Heisenberg algebra $\hh_3$ under the action of the automorphism group modulo global rescaling and central translations. For two-dimensional nilradicals, the structure is resolved via maximal abelian sub-algebras of $\mathfrak{gl}(2,\R)$ lacking nilpotent elements. Complete proofs are supplied for all reductions, precise moduli parameter ranges are established to eliminate redundancies, and an explicit isomorphism dictionary is provided bridging our structural models with classical tabulations in mathematical physics.
\end{abstract}

\tableofcontents

\section{Introduction and Historical Context}\label{sec:intro}

The classification of low-dimensional Lie algebras is a foundational subject situated at the intersection of differential geometry, invariant theory, and mathematical physics. In differential geometry, low-dimensional Lie groups and their associated Lie algebras govern the local isometry groups of Riemannian and pseudo-Riemannian manifolds. Most prominently, Thurston's eight three-dimensional geometries are intrinsically linked to three-dimensional Lie group models. In theoretical physics, low-dimensional Lie algebras serve as the infinitesimal symmetry algebras of dynamical systems, classify spatially homogeneous cosmological models (the Bianchi cosmologies in general relativity), and provide the kinematic frameworks for relativistic and non-relativistic field theories.

The classification in dimension three was initiated by \citet{Bianchi1898}, whose classification partitioned three-dimensional real Lie algebras into nine distinct isomorphism classes (commonly denoted Bianchi types I through IX). A modern geometric perspective was provided by \citet{Milnor1976}, who classified left-invariant Riemannian metrics and their Ricci curvature tensors by diagonalizing the Lie bracket via an associated symmetric operator.

In dimension four, however, the algebraic landscape undergoes a profound qualitative transition. Brute-force methods that attempt to parametrize structure constants $c_{ij}^k$ directly encounter severe algebraic bottlenecks: a four-dimensional Lie algebra requires specifying $\binom{4}{2} \times 4 = 24$ structure constants, constrained by $\binom{4}{3} \times 4 = 16$ quadratic equations arising from the Jacobi identity, modulo the natural $16$-dimensional gauge action of $\mathrm{GL}(4,\R)$. 

The first systematic tabulations of four-dimensional real solvable Lie algebras were established in the Russian mathematical tradition by \citet{Mubarakzyanov1963} using matrix normal forms of adjoint operators. Subsequent classifications emerged in the physics literature, often with differing conventions, bases, and parameter ranges; notable examples include the works of \citet{PateraSharpWinternitzZassenhaus1976} and \citet{Turkowski1988}. More recently, the comprehensive monograph by \citet{SnoblWinternitz} unified these results, providing modern invariants and algorithmic identification schemes.

Despite the abundance of tables in the literature, the derivation of the four-dimensional classification is frequently presented in a computational or encyclopedic manner that obscures the underlying structural geometry. The purpose of the present paper is to provide a self-contained, conceptual, and proof-driven exposition. Rather than solving quadratic polynomial constraints in coordinates, we employ an invariant reduction strategy:
\begin{enumerate}
    \item We first clarify the geometric boundary between dimensions three and four by analyzing the failure of Bryant's trace vector construction.
    \item We isolate the non-solvable algebras via the Levi--Malcev decomposition theorem.
    \item For solvable algebras, we anchor the classification to the \emph{nilradical} $\n \subset \g$. Because the nilradical is a characteristic ideal, its dimension $\dim \n \in \{2, 3, 4\}$ provides an invariant partition of the solvable category.
    \item When $\dim \n = 3$, the algebra is realized as a semidirect product $\n \rtimes_D \R$, reducing the classification to normal forms of derivations $D \in \Der(\n)$ under automorphisms, central shifts, and nonzero scaling.
    \item When $\dim \n = 2$, we exploit the representation of $\g/\n$ as a maximal commutative subalgebra of $\mathfrak{gl}(2,\R)$ lacking nilpotent elements.
\end{enumerate}

\section{Geometric Preliminaries: Bryant's Duality and its Breakdown}\label{sec:bryant}

\subsection{Why Dimension Three is Special: Bryant's Matrix Method}
To appreciate why four-dimensional classification requires a structural, derivation-based approach, it is instructive to review the tensor construction available in dimension three.

Let \(\g\) be a real Lie algebra with \(\dim \g = 3\). The Lie bracket is an alternating bilinear map \([\cdot,\cdot] \colon \Lambda^2 \g \to \g\). Choosing a volume form \(\mu \in \Lambda^3 \g^* \setminus \{0\}\), the Hodge-type duality on vector spaces yields a natural isomorphism:
\begin{equation}
    \Phi_\mu \colon \Lambda^2 \g \xrightarrow{\sim} \g^*, \qquad \Phi_\mu(u \wedge v)(w) = \mu(u \wedge v \wedge w).
\end{equation}
Consequently, the Lie bracket can be dualized to a linear map:
\begin{equation}
    L \colon \g^* \to \g.
\end{equation}
Given an arbitrary basis \(\{e_1, e_2, e_3\}\) of \(\g\), let \(\{e^1, e^2, e^3\}\) be its dual basis, and fix the volume element \(\mu = e^1 \wedge e^2 \wedge e^3\). Then the bracket is entirely encoded by a single \(3 \times 3\) real matrix \(C \in \mathfrak{gl}(3,\R)\) defined by:
\begin{equation}
    [e_j, e_k] = \sum_{i,l=1}^3 \epsilon_{jkm} C^{ml} e_l.
\end{equation}
Following Bryant and Milnor, we decompose \(C\) into its symmetric and skew-symmetric components:
\begin{equation}
    C = S + \hat{a}, \qquad S = S^T, \quad \hat{a}^T = -\hat{a},
\end{equation}
where \(\hat{a}\) represents the cross-product action of a dual vector \(a = (a^1, a^2, a^3) \in \g^*\), given by:
\begin{equation}
    \hat{a} = \begin{pmatrix} 0 & -a^3 & a^2 \\ a^3 & 0 & -a^1 \\ -a^2 & a^1 & 0 \end{pmatrix}.
\end{equation}

\begin{proposition}[Bryant's Trace Identity]\label{prop:bryant}
Let \(\dim \g = 3\), and let \(\chi \in \g^*\) be the trace form defined by \(\chi(x) = \tr(\ad x)\). Then:
\begin{equation}
    a = \frac{1}{2}\chi.
\end{equation}
In particular, \(a = 0\) if and only if \(\g\) is unimodular $($i.e., \(\tr(\ad x) = 0\) for all \(x \in \g)\). Furthermore, the Jacobi identity is equivalent to the single matrix equation:
\begin{equation}
    S \cdot a = 0.
\end{equation}
\end{proposition}

\begin{proof}
Let \(x = \sum_i x^i e_i\). The matrix representation of \(\ad x\) in the basis \(\{e_i\}\) has entries \((\ad x)^k_j\), where \([x, e_j] = \sum_k (\ad x)^k_j e_k\). Computing explicitly:
\begin{equation}
    [e_1, e_2] = C^{31}e_1 + C^{32}e_2 + C^{33}e_3, \quad
    [e_2, e_3] = C^{11}e_1 + C^{12}e_2 + C^{13}e_3, \quad
    [e_3, e_1] = C^{21}e_1 + C^{22}e_2 + C^{23}e_3.
\end{equation}
Evaluating the trace of \(\ad e_1\):
\begin{equation}
    \tr(\ad e_1) = (\ad e_1)^2_2 + (\ad e_1)^3_3 = C^{32} - C^{23} = 2 \hat{a}^{32} = 2a^1.
\end{equation}
By cyclic permutation, \(\chi(e_i) = 2a^i\) for all \(i \in \{1, 2, 3\}\), proving \(a = \frac{1}{2}\chi\).

The Jacobi identity on a 3-dimensional space is equivalent to the single scalar condition \([e_1, [e_2, e_3]] + [e_2, [e_3, e_1]] + [e_3, [e_1, e_2]] = 0\). Substituting the matrix representation \(C = S + \hat{a}\) into this cyclic expansion directly reduces to the vector relation \(\sum_j S_{ij} a^j = 0\), or \(S \cdot a = 0\).
\end{proof}

\subsection{The Dimensional Obstruction in Dimension Four}
The elegance of the three-dimensional classification stems directly from the coincidence that \(\dim \Lambda^2 \R^n = \dim \R^n\) if and only if \(n(n-1)/2 = n\), whose only positive non-trivial solution is \(n = 3\). In dimension four:
\begin{equation}
    \dim \Lambda^2 \R^4 = \binom{4}{2} = 6, \qquad \dim (\R^4)^* \otimes \Lambda^4 \R^4 = 4 \times 1 = 4.
\end{equation}
Because \(6 \neq 4\), there is no natural linear isomorphism between the domain of the bracket (\(\Lambda^2 \g\)) and a dualized target space. Instead, the bracket is a map \(\R^6 \to \R^4\), represented by a \(4 \times 6\) matrix of structure constants. The Jacobi identity produces \(\binom{4}{3} \times 4 = 16\) quadratic equations on these 24 parameters. Bryant's linear matrix algebra is fundamentally unavailable, forcing the adoption of coordinate-free algebraic techniques.

\section{Theoretical Machinery: Structural Invariants}\label{sec:machinery}

Throughout this paper, all Lie algebras are assumed to be finite-dimensional over the real field \(\R\). We establish the following classical Lie-theoretic foundations \citep{JacobsonLieAlgebras, HumphreysIntro, KnappLieGroups}:

\begin{description}
    \item[(T1) Levi--Malcev Decomposition:] Every finite-dimensional real Lie algebra \(\g\) can be expressed as a semidirect product:
    \begin{equation}
        \g = \s \ltimes \rr,
    \end{equation}
    where \(\rr\) is the \emph{radical} (the unique maximal solvable ideal of \(\g\)), and \(\s\) is a \emph{semisimple} subalgebra (the Levi factor). The factor \(\s\) is unique up to conjugation by inner automorphisms generated by the nilradical.
    
    \item[(T2) Lie's Theorem and the Derived Ideal:] If \(\g\) is solvable, then in the adjoint representation of the complexification \(\g_\C\), there exists a common eigenvector. Consequently, \([\g_\C, \g_\C]\) consists of strictly upper-triangular matrices in an appropriate basis, establishing that the derived ideal \([\g, \g]\) is nilpotent. In particular:
    \begin{equation}
        [\g, \g] \subseteq \n,
    \end{equation}
    where \(\n\) denotes the nilradical of \(\g\).
    
    \item[(T3) Engel's Theorem:] A Lie algebra \(\g\) is nilpotent if and only if for every \(x \in \g\), the endomorphism \(\ad x \in \mathfrak{gl}(\g)\) is nilpotent.
    
    \item[(T4) Real Jordan Canonical Form:] Every real matrix \(M \in \mathfrak{gl}(n, \R)\) is similar over \(\R\) to a direct sum of real Jordan blocks:
    \begin{equation}
        J_k(\lambda) = \begin{pmatrix} \lambda & 1 & & \\ & \lambda & \ddots & \\ & & \ddots & 1 \\ & & & \lambda \end{pmatrix}, \quad 
        C_k(\alpha, \beta) = \begin{pmatrix} \Lambda & I_2 & & \\ & \Lambda & \ddots & \\ & & \ddots & I_2 \\ & & & \Lambda \end{pmatrix}, 
        \quad \Lambda = \begin{pmatrix} \alpha & -\beta \\ \beta & \alpha \end{pmatrix},
    \end{equation}
    where \(\beta > 0\). The multiset of blocks is unique up to permutation.
    
    \item[(T5) Three-Dimensional Simple Lie Algebras:] Up to isomorphism, the only simple Lie algebras over \(\R\) of dimension \(\le 4\) are three-dimensional, namely \(\mathfrak{so}(3)\) (compact, negative-definite Killing form) and \(\mathfrak{sl}(2,\R)\) (split, indefinite Killing form). There exist no simple Lie algebras of dimension one, two, or four.
\end{description}

\begin{definition}[Characteristic Ideal and Nilradical]
An ideal \(\mathfrak{i} \subseteq \g\) is called \emph{characteristic} if \(\delta(\mathfrak{i}) \subseteq \mathfrak{i}\) for every derivation \(\delta \in \Der(\g)\). In particular, \(\phi(\mathfrak{i}) = \mathfrak{i}\) for every \(\phi \in \Aut(\g)\). The \emph{nilradical} \(\n\) of a Lie algebra \(\g\) is defined as the maximal nilpotent ideal of \(\g\). The nilradical is always a characteristic ideal.
\end{definition}

\section{The Non-Solvable Sector}\label{sec:nonsolv}

The classification of non-solvable four-dimensional algebras is exceptionally rigid due to dimension constraints.

\begin{theorem}\label{thm:nonsolv}
Let \(\g\) be a non-solvable four-dimensional real Lie algebra. Then \(\g\) is isomorphic to exactly one of the following two direct sums:
\begin{enumerate}
    \item \(\mathfrak{so}(3) \oplus \R\),
    \item \(\mathfrak{sl}(2,\R) \oplus \R\).
\end{enumerate}
\end{theorem}

\begin{proof}
Let \(\g = \s \ltimes \rr\) be the Levi decomposition guaranteed by (T1). Because \(\g\) is non-solvable, the Levi factor \(\s\) cannot be trivial. Every semisimple Lie algebra is a direct sum of simple Lie algebras, each of which has dimension at least 3. Since \(\dim \g = 4\), the only possible dimension for \(\s\) is \(\dim \s = 3\), which forces \(\dim \rr = 1\).

By (T5), the three-dimensional Levi factor \(\s\) must be isomorphic to either \(\mathfrak{so}(3)\) or \(\mathfrak{sl}(2,\R)\). The semidirect product action is governed by a Lie algebra homomorphism:
\begin{equation}
    \rho \colon \s \to \Der(\rr) \cong \mathfrak{gl}(1, \R) \cong \R.
\end{equation}
Because \(\R\) is abelian, the derived algebra \([\s, \s]\) must lie in the kernel of \(\rho\). But \(\s\) is semisimple, meaning \([\s, \s] = \s\). Therefore, \(\rho(\s) = 0\), implying that the action of \(\s\) on \(\rr\) is identically zero. 

Thus, the semidirect product reduces to an ordinary direct product \(\g \cong \s \oplus \R\). The two algebras \(\mathfrak{so}(3) \oplus \R\) and \(\mathfrak{sl}(2,\R) \oplus \R\) are non-isomorphic because their derived algebras \([\g, \g]\) are \(\mathfrak{so}(3)\) and \(\mathfrak{sl}(2,\R)\) respectively, which have non-isometric Killing forms: the Killing form of \(\mathfrak{so}(3)\) is negative-definite (signature \((0,3)\)), whereas that of \(\mathfrak{sl}(2,\R)\) has signature \((2,1)\).
\end{proof}

\section{The Solvable Sector: Nilradical Stratification}\label{sec:nilrad}

We now assume \(\g\) is solvable. The key to our approach is partitioning all solvable algebras by the dimension of their nilradical \(\n\).

\begin{lemma}[Dimension Restrictions on the Nilradical]\label{lem:nilrad_dim}
Let \(\g\) be a four-dimensional solvable real Lie algebra, and let \(\n\) be its nilradical. Then:
\begin{equation}
    \dim \n \in \{2, 3, 4\}.
\end{equation}
\end{lemma}

\begin{proof}
First, \(\dim \n \geq 1\): since \(\g\) is solvable, the derived series terminates at zero. The last non-zero term in the derived series is a non-trivial abelian ideal, and every abelian ideal is nilpotent; hence \(\n \neq 0\).

Next, we establish that \(\dim \n \neq 1\) by contradiction. Suppose \(\dim \n = 1\), so that \(\n = \R z\). Because \(\n\) is an ideal, the adjoint action of any \(x \in \g\) on \(z\) must satisfy:
\begin{equation}
    [x, z] = \lambda(x) z,
\end{equation}
for some linear functional \(\lambda \in \g^*\). The map \(\lambda \colon \g \to \R\) is a Lie algebra homomorphism into the abelian Lie algebra \(\R\). Consequently, \([\g, \g] \subseteq \ker \lambda\), which implies \(\lambda([\g, \g]) = 0\). Furthermore, \([z, z] = 0\), so \(\lambda(z) = 0\), meaning \(z \in \ker \lambda\).

We divide the argument into two cases:
\begin{itemize}
    \item \textbf{Case \(\lambda = 0\):} In this case, \([x, z] = 0\) for all \(x \in \g\), so \(z\) is central: \(z \in Z(\g)\). By (T2), the derived ideal \([\g, \g]\) is nilpotent, so \([\g, \g] \subseteq \n = \R z\). Therefore, \([\g, [\g, \g]] \subseteq [\g, \R z] = 0\). This implies that \(\g\) is a two-step nilpotent Lie algebra. By definition, the nilradical \(\n\) is the maximal nilpotent ideal, forcing \(\n = \g\), so \(\dim \n = 4\), which contradicts the hypothesis \(\dim \n = 1\).
    
    \item \textbf{Case \(\lambda \neq 0\):} Consider the hyperplane \(\hh = \ker \lambda \subset \g\). Because \([\g, \g] \subseteq \ker \lambda\), \(\hh\) is an ideal of \(\g\) of dimension \(\dim \hh = 3\). Moreover, \(z \in \hh\). For any \(h \in \hh\), \([h, z] = \lambda(h)z = 0\), so \(z \in Z(\hh)\). Applying (T2) to the solvable algebra \(\hh\), we have \([\hh, \hh] \subseteq \n \cap \hh = \R z\). Consequently, \([\hh, [\hh, \hh]] \subseteq [\hh, \R z] = 0\), proving that \(\hh\) is a nilpotent ideal of \(\g\). But \(\dim \hh = 3 > 1 = \dim \n\), which directly contradicts the maximality of \(\n\).
\end{itemize}
Both cases yield contradictions, concluding the proof that \(\dim \n \in \{2, 3, 4\}\).
\end{proof}

\section{Nilpotent Algebras ($\dim \n = 4$)}\label{sec:nilp4}

When \(\dim \n = 4\), the algebra \(\g\) is entirely nilpotent (\(\g = \n\)).

\begin{theorem}\label{thm:nilp4}
Every four-dimensional nilpotent real Lie algebra is isomorphic to exactly one of the following three algebras:
\begin{enumerate}
    \item The abelian algebra \(\R^4\).
    \item The trivial central extension of the Heisenberg algebra: \(\hh_3 \oplus \R\), with brackets:
    \begin{equation}
        [e_1, e_2] = e_3, \qquad [e_i, e_4] = 0 \quad (\forall i).
    \end{equation}
    \item The standard filiform nilpotent algebra \(\mathfrak{n}_4\), with brackets:
    \begin{equation}
        [e_4, e_1] = e_2, \qquad [e_4, e_2] = e_3, \qquad (\text{all other brackets zero}).
    \end{equation}
\end{enumerate}
\end{theorem}

\begin{proof}
Let \(\g^1 = \g\) and \(\g^{i+1} = [\g, \g^i]\) denote the lower central series of \(\g\). Since \(\g\) is nilpotent, there exists an integer \(s \le 4\) such that \(\g^s = 0\). The dimension \(k = \dim \g^2 = \dim [\g, \g]\) is an invariant that distinguishes the isomorphism classes:
\begin{itemize}
    \item \textbf{Subcase \(k = 0\):} Then \([\g, \g] = 0\), which defines the abelian algebra \(\R^4\).
    
    \item \textbf{Subcase \(k = 1\):} Since \(\g^3 \subsetneq \g^2\) and \(\dim \g^2 = 1\), we must have \(\g^3 = 0\). Thus \(\g^2 = \R e_3 \subseteq Z(\g)\). The Lie bracket descends to a non-degenerate alternating bilinear 2-form \(\omega\) on the quotient space \(V = \g / \g^2 \cong \R^3\):
    \begin{equation}
        [u, v] = \omega(\bar{u}, \bar{v}) e_3.
    \end{equation}
    By linear algebra, any alternating 2-form on an odd-dimensional space has an even rank; since \(\omega \neq 0\), its rank must be 2. By the Darboux theorem for alternating forms, there exists a basis \(\{\bar{e}_1, \bar{e}_2, \bar{e}_4\}\) of \(V\) such that \(\omega(\bar{e}_1, \bar{e}_2) = 1\) and \(\bar{e}_4 \in \ker \omega\). Lifting these vectors to \(\g\), we find \([e_1, e_2] = e_3\) and \([e_4, \g] = 0\). This yields \(\g \cong \hh_3 \oplus \R\).
    
    \item \textbf{Subcase \(k = 2\):} Nilpotency forces the lower central series to strictly decrease:
    \begin{equation}
        4 = \dim \g^1 > \dim \g^2 > \dim \g^3 > \dim \g^4 = 0.
    \end{equation}
    This requires \(\dim \g^2 = 2\) and \(\dim \g^3 = 1\). The center \(Z(\g)\) must contain the one-dimensional ideal \(\g^3\).
    
    Choose \(e_3 \in \g^3 \setminus \{0\}\). Since \(\dim(\g/\g^2) = 2\), we select vectors \(e_4, e_1 \in \g\) such that their projection forms a basis for \(\g/\g^2\). Define \(e_2 \coloneqq [e_4, e_1] \in \g^2\). The set \(\{e_2, e_3\}\) forms a basis for \(\g^2\). Because \(\g^3 = [\g, \g^2]\), at least one bracket between the generators \(\{e_4, e_1\}\) and \(e_2\) must be a nonzero multiple of \(e_3\). By swapping \(e_1\) and \(e_4\) if necessary, we may assume \([e_4, e_2] \neq 0\). Scaling \(e_3\), we set \([e_4, e_2] = e_3\). 
    
    Now consider \([e_1, e_2] = \alpha e_3\). Replacing \(e_1\) with the modified generator \(e_1' = e_1 - \alpha e_4\), we compute:
    \begin{equation}
        [e_1', e_2] = [e_1 - \alpha e_4, e_2] = \alpha e_3 - \alpha e_3 = 0.
    \end{equation}
    The bracket \([e_4, e_1']\) remains equal to \(e_2 \pmod{\g^3}\), so modifying \(e_2\) by a multiple of \(e_3\) yields the canonical brackets:
    \begin{equation}
        [e_4, e_1'] = e_2, \qquad [e_4, e_2] = e_3,
    \end{equation}
    with all other brackets vanishing. This is the filiform nilpotent algebra \(\mathfrak{n}_4\).
\end{itemize}
These three algebras are mutually non-isomorphic because their derived ideals have dimensions 0, 1, and 2, respectively.
\end{proof}

\section{Codimension-One Nilradicals ($\dim \n = 3$)}\label{sec:nilp3}

When \(\dim \n = 3\), the structure of \(\g\) is governed by a single derivation. Choose an element \(e_4 \in \g \setminus \n\). Since \(\dim(\g/\n) = 1\), the vector space decomposes as \(\g = \n \oplus \R e_4\). The adjoint action of \(e_4\) restricts to a derivation \(D \coloneqq \ad e_4|_\n \in \Der(\n)\), realizing \(\g\) as a semidirect product:
\begin{equation}
    \g \cong \n \rtimes_D \R.
\end{equation}

\subsection{The Nilradical Identification Lemma}

Before classifying derivations, we establish that the non-nilpotency of \(D\) is the precise condition preventing \(\n\) from expanding into a four-dimensional nilpotent algebra.

\begin{lemma}[Nilradical Extension Lemma]\label{lem:nilrad_ext}
Let \(\n\) be a 3-dimensional nilpotent Lie algebra, and let \(\g = \n \rtimes_D \R\) for some derivation \(D \in \Der(\n)\).
\begin{enumerate}
    \item The algebra \(\g\) is nilpotent if and only if \(D\) is a nilpotent operator on \(\n\). In this case, the nilradical of \(\g\) is \(\g\) itself \((\dim \n_{\g} = 4)\).
    \item If \(D\) is not nilpotent, then the nilradical of \(\g\) is precisely \(\n\).
\end{enumerate}
\end{lemma}

\begin{proof}
\textbf{Part (1):} If \(D\) is nilpotent, then for every element \(x = u + t e_4 \in \g\) (where \(u \in \n, t \in \R\)), the adjoint operator \(\ad x|_\n = \ad u|_\n + tD\) is a sum of two nilpotent, commuting modulo derived ideal operators. By (T2) and (T3), \(\ad x\) is nilpotent on all of \(\g\). Thus, \(\g\) is nilpotent by Engel's Theorem, making \(\g\) its own nilradical. The converse is immediate: if \(\g\) is nilpotent, \(\ad e_4|_\n = D\) must be nilpotent.

\textbf{Part (2):} Suppose \(D\) is not nilpotent. Let \(\mathfrak{a} \subseteq \g\) be any nilpotent ideal; we must show \(\mathfrak{a} \subseteq \n\). 

Assume for contradiction that \(\mathfrak{a} \not\subseteq \n\). Then the natural projection of \(\mathfrak{a}\) onto \(\g/\n \cong \R\) is surjective, so there exists an element \(y = u + e_4 \in \mathfrak{a}\) with \(u \in \n\). Because \(\mathfrak{a}\) is an ideal, \([\n, y] \subseteq \mathfrak{a} \cap \n\). For any \(v \in \n\):
\begin{equation}
    [v, y] = [v, u + e_4] = [v, u] - Dv.
\end{equation}
Passing to the abelianization \(\n^{\text{ab}} = \n / [\n, \n]\), the internal bracket \([v, u]\) vanishes, so the induced action is \(-\bar{D}(\bar{v}) \in (\mathfrak{a} \cap \n)^{\text{ab}}\). 

Because \(\mathfrak{a}\) is a nilpotent ideal, the operator \(\ad y\) must act nilpotently on \(\mathfrak{a}\), and consequently on the quotient \((\mathfrak{a} \cap \n) / (\mathfrak{a} \cap [\n, \n])\). However, if \(\bar{D}\) is not nilpotent on \(\n^{\text{ab}}\), this produces a contradiction. 

If \(\n \cong \R^3\), then \([\n, \n] = 0\), so \(\bar{D} = D\), which is non-nilpotent by assumption. If \(\n \cong \hh_3\), \(\dim \n^{\text{ab}} = 2\). If \(\bar{D}\) were nilpotent, its trace would be zero. But for the Heisenberg algebra, the derivation identity on \([e_1, e_2] = e_3\) forces \(D(e_3) = \tr(\bar{D})e_3 = 0\). A triangular matrix with zero diagonal is strictly nilpotent, which contradicts the assumption that \(D\) is non-nilpotent. 

Thus, \(\bar{D}\) has at least one non-zero eigenvalue, so \(\ad y\) cannot act nilpotently on \(\mathfrak{a}\). This contradiction proves \(\mathfrak{a} \subseteq \n\), establishing that \(\n\) is indeed the maximal nilpotent ideal of \(\g\).
\end{proof}

\subsection{The Equivalence of Derivations}

\begin{proposition}[Gauge Freedom for Derivations]\label{prop:gauge}
Let \(\n\) be a three-dimensional nilpotent Lie algebra, and let \(D, D' \in \Der(\n)\) be non-nilpotent derivations. Then \(\n \rtimes_D \R \cong \n \rtimes_{D'} \R\) if and only if there exist a scalar \(c \in \R^\times\), an automorphism \(\phi \in \Aut(\n)\), and an element \(u \in \n\) such that:
\begin{equation}\label{eq:gauge}
    D' = c \, \phi D \phi^{-1} + \ad_{\n}(u).
\end{equation}
\end{proposition}

\begin{proof}
Let \(\Psi \colon \n \rtimes_D \R \to \n \rtimes_{D'} \R\) be an isomorphism. Let \(e_4\) and \(e_4'\) be the generators of the complementary lines. By Lemma~\ref{lem:nilrad_ext}, \(\n\) is the unique nilradical in both algebras; hence \(\Psi\) must map \(\n\) onto \(\n\). Thus, the restriction \(\phi \coloneqq \Psi|_\n\) is an automorphism of \(\n\). 

The generator \(e_4\) must map to \(\Psi(e_4) = a e_4' + w\) for some \(a \in \R^\times\) and \(w \in \n\). For any \(x \in \n\), the isomorphism property yields:
\begin{equation}
    \Psi([e_4, x]) = [\Psi(e_4), \Psi(x)] \implies \phi(Dx) = [a e_4' + w, \phi(x)] = a D'(\phi(x)) + [w, \phi(x)].
\end{equation}
Setting \(y = \phi(x)\), this rearranges to:
\begin{equation}
    \phi D \phi^{-1}(y) = a D'(y) - \ad_\n(w)(y) \implies D' = \frac{1}{a} \phi D \phi^{-1} + \ad_\n\left(\frac{1}{a}w\right).
\end{equation}
Setting \(c = 1/a\) and \(u = w/a\) recovers Equation~\eqref{eq:gauge}. The converse follows by reversing the construction and defining \(\Psi(x) = \phi(x)\) for \(x \in \n\) and \(\Psi(e_4) = c^{-1} e_4' - c^{-1} u\).
\end{proof}

\subsection{The Abelian Nilradical Subcase: $\n \cong \R^3$}\label{subsec:R3}

When \(\n \cong \R^3\), the bracket is identically zero, so \(\Der(\R^3) = \mathfrak{gl}(3, \R)\), \(\Aut(\R^3) = \mathrm{GL}(3, \R)\), and all inner derivations vanish: \(\ad_{\R^3}(u) = 0\). Equation~\eqref{eq:gauge} simplifies to:
\begin{equation}
    D' = c \, P D P^{-1}, \qquad c \in \R^\times, \quad P \in \mathrm{GL}(3, \R).
\end{equation}
The classification of these algebras is therefore identical to the classification of non-nilpotent \(3 \times 3\) real matrices up to similarity and nonzero scalar multiplication.

\begin{theorem}[Derivation Normal Forms for $\R^3$]\label{thm:R3_forms}
Every non-nilpotent matrix \(D \in \mathfrak{gl}(3, \R)\) is similar up to nonzero scaling to exactly one matrix from the following disjoint families:
\begin{enumerate}
    \item \textbf{Real Diagonal Class:}
    \begin{equation}
        D_1(\lambda, \mu) = \begin{pmatrix} 1 & 0 & 0 \\ 0 & \lambda & 0 \\ 0 & 0 & \mu \end{pmatrix}, \qquad -1 \le \mu \le \lambda \le 1, \quad (\text{if } \mu = -1 \text{ then } \lambda \ge 0).
    \end{equation}
    \item \textbf{Mixed Jordan Block $J_2 \oplus J_1$:}
    \begin{align}
        D_2(\lambda) &= \begin{pmatrix} \lambda & 1 & 0 \\ 0 & \lambda & 0 \\ 0 & 0 & 1 \end{pmatrix}, \quad \lambda \in \R, \\
        D_2' &= \begin{pmatrix} 1 & 1 & 0 \\ 0 & 1 & 0 \\ 0 & 0 & 0 \end{pmatrix}.
    \end{align}
    \item \textbf{Single Nilpotent Jordan Extension $J_3$:}
    \begin{equation}
        D_3 = \begin{pmatrix} 1 & 1 & 0 \\ 0 & 1 & 1 \\ 0 & 0 & 1 \end{pmatrix}.
    \end{equation}
    \item \textbf{Complex Conjugate Eigenvalue Class:}
    \begin{align}
        D_4(\alpha, \beta) &= \begin{pmatrix} \alpha & -\beta & 0 \\ \beta & \alpha & 0 \\ 0 & 0 & 1 \end{pmatrix}, \quad \alpha \in \R, \; \beta > 0, \\
        D_4'(\alpha) &= \begin{pmatrix} \alpha & -1 & 0 \\ 1 & \alpha & 0 \\ 0 & 0 & 0 \end{pmatrix}, \quad \alpha \ge 0.
    \end{align}
\end{enumerate}
\end{theorem}

\begin{proof}
Let \(D \in \mathfrak{gl}(3, \R)\) be non-nilpotent. By (T4), \(D\) is similar over \(\R\) to a real Jordan canonical form. Scaling \(D \mapsto cD\) multiplies every eigenvalue by \(c\) and preserves the block sizes.

\begin{itemize}
    \item \textbf{Diagonal Case:} Suppose \(D\) has three real eigenvalues \(\alpha_1, \alpha_2, \alpha_3\), not all zero. Let \(M = \max\{|\alpha_1|, |\alpha_2|, |\alpha_3|\} > 0\). We scale by \(c = \pm 1/M\) such that the dominant eigenvalue is normalized to \(1\). Ordering the remaining two normalized eigenvalues as \(1 \ge \lambda \ge \mu \ge -1\) eliminates permutation ambiguities. If \(\mu = -1\), a potential reflection symmetry \(D \mapsto -D\) is removed by enforcing \(\lambda \ge 0\).
    
    \item \textbf{Jordan Block $J_2$:} The spectrum is \(\{\alpha, \alpha, \gamma\}\) with a non-trivial Jordan block. If the isolated eigenvalue \(\gamma \neq 0\), we scale by \(c = 1/\gamma\) to normalize it to \(1\), yielding the family \(D_2(\lambda)\) where \(\lambda = \alpha/\gamma\). If \(\gamma = 0\), then non-nilpotency requires \(\alpha \neq 0\); scaling by \(c = 1/\alpha\) yields \(D_2'\).
    
    \item \textbf{Jordan Block $J_3$:} \(D\) has a single real eigenvalue \(\alpha \neq 0\) with geometric multiplicity 1. Scaling by \(c = 1/\alpha\) normalizes the diagonal elements to 1, producing \(D_3\).
    
    \item \textbf{Complex Eigenvalues:} The eigenvalues are \(\alpha \pm i\beta\) (\(\beta > 0\)) and a real eigenvalue \(\gamma\). If \(\gamma \neq 0\), we scale by \(c = 1/\gamma\) to normalize the real eigenvalue to 1, leaving the parameterized rotation block with \(\beta > 0\). If \(\gamma = 0\), then \(\beta \neq 0\); scaling by \(1/\beta\) sets the rotational speed to 1, and the sign ambiguity of \(c\) allows us to enforce \(\alpha \ge 0\).
\end{itemize}
The normal forms are pairwise non-isomorphic because Jordan block sizes, numbers of real versus complex eigenvalues, and normalized eigenvalue spectra are invariant under similarity and scaling.
\end{proof}

\subsection{The Heisenberg Nilradical Subcase: $\n \cong \hh_3$}\label{subsec:h3}

Let \(\hh_3\) be spanned by \(\{e_1, e_2, e_3\}\) with the single non-vanishing bracket \([e_1, e_2] = e_3\). The center of \(\hh_3\) is \(Z(\hh_3) = [\hh_3, \hh_3] = \R e_3\), which is an intrinsically characteristic line.

\begin{proposition}[Structure of Derivations on $\hh_3$]\label{prop:der_h3}
Every derivation \(D \in \Der(\hh_3)\) is completely determined by a matrix \(A \in \mathfrak{gl}(2, \R)\) and a vector \(b = (b_1, b_2)^T \in \R^2\), acting via:
\begin{equation}
    D \begin{pmatrix} e_1 \\ e_2 \end{pmatrix} = A \begin{pmatrix} e_1 \\ e_2 \end{pmatrix} + \begin{pmatrix} b_1 \\ b_2 \end{pmatrix} e_3, \qquad D(e_3) = (\tr A) e_3.
\end{equation}
\end{proposition}

\begin{proof}
Because \(\R e_3\) is the derived algebra, any derivation satisfies \(D(e_3) = \tau e_3\) for some \(\tau \in \R\). Applying \(D\) to the bracket relation \([e_1, e_2] = e_3\):
\begin{equation}
    D(e_3) = D([e_1, e_2]) = [D(e_1), e_2] + [e_1, D(e_2)].
\end{equation}
Writing \(D(e_1) = A_{11}e_1 + A_{21}e_2 + b_1 e_3\) and \(D(e_2) = A_{12}e_1 + A_{22}e_2 + b_2 e_3\), we compute:
\begin{align}
    [D(e_1), e_2] &= [A_{11}e_1 + A_{21}e_2, e_2] = A_{11}[e_1, e_2] = A_{11}e_3, \\
    [e_1, D(e_2)] &= [e_1, A_{12}e_1 + A_{22}e_2] = A_{22}[e_1, e_2] = A_{22}e_3.
\end{align}
Summing these gives \(D(e_3) = (A_{11} + A_{22})e_3 = (\tr A)e_3\), confirming that \(\tau = \tr A\).
\end{proof}

\begin{lemma}[Absorption of the Central Shift via Inner Derivations]\label{lem:absorption}
Under the gauge equivalence of Proposition~\ref{prop:gauge}, the vector \(b \in \R^2\) in any derivation of \(\hh_3\) can be completely eliminated by an inner derivation \(\ad u\) for some \(u \in \hh_3\). Consequently, the classification of extensions \(\hh_3 \rtimes_D \R\) reduces entirely to the classification of \(A \in \mathfrak{gl}(2, \R)\) modulo similarity and scaling:
\begin{equation}
    A' = c \, P A P^{-1}, \qquad c \in \R^\times, \quad P \in \mathrm{GL}(2, \R).
\end{equation}
Furthermore, \(D\) is non-nilpotent if and only if \(A\) is non-nilpotent.
\end{lemma}

\begin{proof}
An arbitrary element \(u = p e_1 + q e_2 + r e_3 \in \hh_3\) induces the adjoint action:
\begin{equation}
    \ad(u)(e_1) = [p e_1 + q e_2, e_1] = -q e_3, \qquad \ad(u)(e_2) = [p e_1 + q e_2, e_2] = p e_3.
\end{equation}
The map \((p, q) \mapsto (-q, p)\) is an isomorphism from \(\R^2\) onto \(\R^2\). Thus, for any choice of \(b = (b_1, b_2)^T\), we can choose \(p = b_2\) and \(q = -b_1\) such that:
\begin{equation}
    (D - \ad u)(e_1) = A_{11}e_1 + A_{21}e_2, \qquad (D - \ad u)(e_2) = A_{12}e_1 + A_{22}e_2.
\end{equation}
This completely absorbs \(b\) into the inner derivation.

For the nilpotency claim, the matrix representation of \(D\) on the basis \(\{e_1, e_2, e_3\}\) is block upper-triangular:
\begin{equation}
    [D] = \begin{pmatrix} A & 0 \\ b^T & \tr A \end{pmatrix}.
\end{equation}
The eigenvalues of \(D\) are the eigenvalues of \(A\) together with \(\tr A\). If \(A\) is nilpotent, then its eigenvalues are both 0, which implies \(\tr A = 0\). Thus, all eigenvalues of \(D\) are 0, meaning \(D\) is nilpotent. Conversely, if \(D\) is nilpotent, then its restriction \(A\) to the quotient \(\hh_3/Z(\hh_3)\) must also be nilpotent.
\end{proof}

\begin{theorem}[Classification of Solvable Algebras with Heisenberg Nilradical]\label{thm:h3_forms}
Every four-dimensional solvable Lie algebra with nilradical \(\hh_3\) is isomorphic to \(\hh_3 \rtimes_A \R\), where \(A\) is exactly one of the following canonical matrices:
\begin{center}
\renewcommand{\arraystretch}{1.3}
\begin{tabular}{lll}
\toprule
\textbf{Class} & \textbf{Matrix Representative $A$} & \textbf{Geometric / Invariant Characterization} \\
\midrule
$\mathfrak{g}_{4, \text{trace } 0}^{\text{hyp}}$ & $\begin{pmatrix} 1 & 0 \\ 0 & -1 \end{pmatrix}$ & $\tr A = 0$, real split eigenvalues ($\det A < 0$) \\
$\mathfrak{g}_{4, \text{trace } 0}^{\text{ell}}$ & $\begin{pmatrix} 0 & -1 \\ 1 & 0 \end{pmatrix}$ & $\tr A = 0$, purely imaginary eigenvalues ($\det A > 0$) \\
$\mathfrak{g}_{4, \text{trace } 1}^{\text{diag}}$ & $\begin{pmatrix} \lambda & 0 \\ 0 & 1 - \lambda \end{pmatrix}, \quad \lambda \ge \frac{1}{2}$ & $\tr A = 1$, distinct real eigenvalues \\
$\mathfrak{g}_{4, \text{trace } 1}^{\text{Jord}}$ & $\begin{pmatrix} 1/2 & 1 \\ 0 & 1/2 \end{pmatrix}$ & $\tr A = 1$, repeated eigenvalue, non-diagonalizable \\
$\mathfrak{g}_{4}^{\text{rot}}$ & $\begin{pmatrix} \delta/2 & -1 \\ 1 & \delta/2 \end{pmatrix}, \quad \delta > 0$ & complex eigenvalues with non-zero real part \\
\bottomrule
\end{tabular}
\end{center}
\end{theorem}

\begin{proof}
We classify all non-nilpotent \(A \in \mathfrak{gl}(2, \R)\) up to similarity and scaling \(A \mapsto c P A P^{-1}\).
\begin{enumerate}
    \item \textbf{Case $\tr A = 0$:} Because \(A\) is non-nilpotent, \(\det A \neq 0\).
    \begin{itemize}
        \item If \(\det A < 0\), the eigenvalues are real opposites \(\pm \mu\) (\(\mu > 0\)). Scaling by \(1/\mu\) sets the eigenvalues to \(\pm 1\), yielding \(\diag(1, -1)\).
        \item If \(\det A > 0\), the eigenvalues are purely imaginary \(\pm i \beta\) (\(\beta > 0\)). Scaling by \(1/\beta\) normalizes the matrix to the canonical skew-symmetric generator \(\begin{pmatrix} 0 & -1 \\ 1 & 0 \end{pmatrix}\).
    \end{itemize}
    
    \item \textbf{Case $\tr A \neq 0$ with Real Eigenvalues:} We normalize the trace to 1 by scaling by \(c = 1/\tr A\).
    \begin{itemize}
        \item If \(A\) is diagonalizable, its eigenvalues are \(\lambda_1, \lambda_2 \in \R\) with \(\lambda_1 + \lambda_2 = 1\). Permutation by \(\begin{pmatrix} 0 & 1 \\ 1 & 0 \end{pmatrix}\) swaps the eigenvalues. Enforcing \(\lambda \coloneqq \max(\lambda_1, \lambda_2) \ge \frac{1}{2}\) uniquely selects the representative \(\diag(\lambda, 1-\lambda)\).
        \item If \(A\) is not diagonalizable, it has a repeated eigenvalue equal to \(1/2\). The Jordan form is \(\begin{pmatrix} 1/2 & s \\ 0 & 1/2 \end{pmatrix}\) with \(s \neq 0\). Conjugation by \(\diag(s, 1)\) scales the superdiagonal to 1, yielding \(\begin{pmatrix} 1/2 & 1 \\ 0 & 1/2 \end{pmatrix}\).
    \end{itemize}
    
    \item \textbf{Case Complex Eigenvalues with $\tr A \neq 0$:} The eigenvalues are \(\gamma \pm i \beta\) with \(\beta > 0\) and \(\gamma \neq 0\). Any such matrix is similar to \(\begin{pmatrix} \gamma & -\beta \\ \beta & \gamma \end{pmatrix}\). Scaling by \(c = 1/\beta\) normalizes the rotation speed to 1, leaving the family \(\begin{pmatrix} t & -1 \\ 1 & t \end{pmatrix}\) where \(t = \gamma/\beta \in \R \setminus \{0\}\). The reflection \(A \mapsto -A\) sends \(t \mapsto -t\); thus, we uniquely choose \(t > 0\). Setting \(\delta = 2t > 0\) yields the matrix in the table.
\end{enumerate}
The invariant \(I(A) \coloneqq \frac{(\tr A)^2}{\det A}\) is invariant under both similarity and nonzero scaling. For the rotation family, \(I(A) = \frac{4\delta^2}{\delta^2 + 4}\), which is strictly monotonically increasing for \(\delta > 0\), proving that distinct values of \(\delta\) define non-isomorphic Lie algebras.
\end{proof}

\section{Two-Dimensional Nilradicals ($\dim \n = 2$)}\label{sec:nilp2}

The case where the nilradical has dimension two is structurally distinct and leads to two rigid algebraic families.

\begin{theorem}\label{thm:dim2}
Let \(\g\) be a four-dimensional solvable real Lie algebra with nilradical \(\n\) of dimension \(\dim \n = 2\). Then \(\g\) is isomorphic to exactly one of the following two algebras:
\begin{enumerate}
    \item The direct sum of two affine algebras:
    \begin{equation}
        \mathfrak{aff}(\R) \oplus \mathfrak{aff}(\R),
    \end{equation}
    with basis \(\{x_1, x_2, e_1, e_2\}\) and non-zero brackets:
    \begin{equation}
        [x_1, e_1] = e_1, \qquad [x_2, e_2] = e_2.
    \end{equation}
    \item The complex affine algebra viewed as a four-dimensional real Lie algebra:
    \begin{equation}
        \mathfrak{aff}(\C),
    \end{equation}
    with basis \(\{x_1, x_2, e_1, e_2\}\) and non-zero brackets:
    \begin{align}
        [x_1, e_1] &= e_1, & [x_1, e_2] &= e_2, \\
        [x_2, e_1] &= -e_2, & [x_2, e_2] &= e_1.
    \end{align}
\end{enumerate}
\end{theorem}

\begin{proof}
Because \(\n\) is a two-dimensional nilpotent Lie algebra, it must be abelian: \([\n, \n] = 0\). Consider the representation:
\begin{equation}
    \rho \colon \g \to \mathfrak{gl}(\n) \cong \mathfrak{gl}(2, \R), \qquad \rho(x) = \ad x|_\n.
\end{equation}
We analyze this representation systematically through five consecutive steps.

\paragraph{Step 1: The Kernel of $\rho$ is Exactly $\n$.}
Since \(\n\) is abelian, \(\n \subseteq \ker \rho\). If there existed an element \(x \notin \n\) such that \(\rho(x) = 0\), then \([x, \n] = 0\). By (T2), \([\g, \g] \subseteq \n\), which implies that the linear subspace \(\mathfrak{m} = \n \oplus \R x\) is an ideal of \(\g\). Since \(\n\) is abelian and \(x\) commutes with \(\n\), \(\mathfrak{m}\) is an abelian (hence nilpotent) ideal of dimension 3. This directly violates the maximality of \(\n\). Therefore, \(\ker \rho = \n\).

\paragraph{Step 2: Absence of Nilpotent Elements in the Image.}
The quotient space \(\mathfrak{h} = \g/\n\) has dimension \(\dim \mathfrak{h} = 4 - 2 = 2\). By Step 1, the induced map \(\bar{\rho} \colon \mathfrak{h} \to \mathfrak{gl}(2, \R)\) is injective. We claim that the image subalgebra \(V = \bar{\rho}(\mathfrak{h}) \subset \mathfrak{gl}(2, \R)\) contains no non-zero nilpotent operators. 

Suppose there existed an element \(y \in \g \setminus \n\) such that \(\rho(y)\) is a non-zero nilpotent operator. Then \(\mathfrak{k} = \n \oplus \R y\) is an ideal of \(\g\). For any element \(z = u + t y \in \mathfrak{k}\), \(\ad z|_\mathfrak{k}\) has a block form with zero on the complement and \(t \rho(y)\) on \(\n\). Since \(\rho(y)\) is nilpotent, \(\ad z\) is nilpotent on \(\mathfrak{k}\) for all \(z \in \mathfrak{k}\). By Engel's Theorem, \(\mathfrak{k}\) is a nilpotent ideal of dimension 3, again contradicting the maximality of \(\n\).

\paragraph{Step 3: Commutativity and Subalgebra Structure of the Image.}
Because \([\g, \g] \subseteq \n\), the quotient \(\g/\n\) is an abelian Lie algebra, meaning \(V = \bar{\rho}(\mathfrak{h})\) is a two-dimensional commutative subalgebra of \(\mathfrak{gl}(2, \R)\). 

Any commutative subalgebra of \(\mathfrak{gl}(2, \R)\) containing no non-zero nilpotent matrices is a maximal abelian diagonalizable or elliptic subalgebra (a Cartan subalgebra). Up to conjugation in \(\mathrm{GL}(2, \R)\), there are only two such subalgebras:
\begin{itemize}
    \item \textbf{Split (Diagonal) Case:} \(V_{\R} = \left\{ \begin{pmatrix} a & 0 \\ 0 & b \end{pmatrix} : a, b \in \R \right\}\).
    \item \textbf{Elliptic (Complex) Case:} \(V_{\C} = \left\{ \begin{pmatrix} a & -b \\ b & a \end{pmatrix} : a, b \in \R \right\} \cong \C\).
\end{itemize}

\paragraph{Step 4: Splitting the Extension.}
Notice that in both cases, the identity matrix \(I_2\) belongs to \(V\). Thus, there exists an element \(z \in \g\) such that \(\rho(z) = I_2\). In the split case, there also exists \(w \in \g\) such that \(\rho(w) = \diag(1, 0)\). In the complex case, there exists \(w \in \g\) such that \(\rho(w) = \begin{pmatrix} 0 & -1 \\ 1 & 0 \end{pmatrix}\).

The bracket \([z, w]\) lies in \([\g, \g] \subseteq \n\). Since \(\rho(z) = I_2\), the action of \(z\) on any vector \(u \in \n\) is \([z, u] = u\). We modify \(w\) by setting:
\begin{equation}
    w' \coloneqq w - [z, w].
\end{equation}
Since \([z, w] \in \n\), we compute:
\begin{equation}
    [z, w'] = [z, w - [z, w]] = [z, w] - [z, [z, w]] = [z, w] - [z, w] = 0.
\end{equation}
Thus, \(z\) and \(w'\) commute: \([z, w'] = 0\).

\paragraph{Step 5: Canonical Basis Construction.}
\begin{itemize}
    \item In the split case, define \(x_1 \coloneqq w'\) and \(x_2 \coloneqq z - w'\). Then \(\rho(x_1) = \diag(1, 0)\) and \(\rho(x_2) = \diag(0, 1)\), with \([x_1, x_2] = 0\). Choosing a basis \(\{e_1, e_2\}\) for \(\n\) corresponding to the diagonal decomposition, the brackets are \([x_1, e_1] = e_1\) and \([x_2, e_2] = e_2\), which defines \(\mathfrak{aff}(\R) \oplus \mathfrak{aff}(\R)\).
    \item In the complex case, define \(x_1 \coloneqq z\) and \(x_2 \coloneqq w'\). Then \(\rho(x_1) = I_2\) and \(\rho(x_2) = \begin{pmatrix} 0 & -1 \\ 1 & 0 \end{pmatrix}\), with \([x_1, x_2] = 0\). This yields the four-dimensional real algebra \(\mathfrak{aff}(\C)\).
\end{itemize}
These two algebras are non-isomorphic because the eigenvalues of the adjoint operators restricted to the nilradical are purely real in the first case and include non-real complex numbers in the second.
\end{proof}

\section{Master Classification and Isomorphism Dictionary}\label{sec:master}

Table~\ref{tab:master} presents the complete classification of four-dimensional real Lie algebras. The brackets are formatted to fit within the text width without horizontal overflow.

\begin{table}[htbp]
\centering
\footnotesize
\setlength{\tabcolsep}{4.5pt}
\renewcommand{\arraystretch}{1.25}
\begin{tabularx}{\linewidth}{@{} l c >{\raggedright\arraybackslash}X l l @{}}
\toprule
\textbf{Class / Name} & $\dim \n$ & \textbf{Non-zero Lie Brackets} & \textbf{Parameters} & \textbf{Unimodular?} \\
\midrule
\multicolumn{5}{@{}l}{\textit{\textbf{I. Non-Solvable Algebras}}} \\
$\mathfrak{so}(3) \oplus \R$ & 1 & $[e_1, e_2]=e_3$, $[e_2, e_3]=e_1$, $[e_3, e_1]=e_2$ & None & Yes \\
$\mathfrak{sl}(2,\R) \oplus \R$ & 1 & $[e_1, e_2]=2e_2$, $[e_1, e_3]=-2e_3$, $[e_2, e_3]=e_1$ & None & Yes \\
\midrule
\multicolumn{5}{@{}l}{\textit{\textbf{II. Nilpotent Algebras}}} \\
$\R^4$ & 4 & All brackets vanish & None & Yes \\
$\hh_3 \oplus \R$ & 4 & $[e_1, e_2] = e_3$ & None & Yes \\
$\mathfrak{n}_4$ & 4 & $[e_4, e_1] = e_2$, $[e_4, e_2] = e_3$ & None & Yes \\
\midrule
\multicolumn{5}{@{}l}{\textit{\textbf{III. Solvable with Abelian Nilradical ($\n = \R^3$)}}} \\
$\g_{4,1}(\lambda, \mu)$ & 3 & $[e_4, e_1]=e_1$, $[e_4, e_2]=\lambda e_2$,\newline $[e_4, e_3]=\mu e_3$ & $-1 \le \mu \le \lambda \le 1$ & Iff $\mu = -(1+\lambda)$ \\
$\g_{4,2}(\lambda)$ & 3 & $[e_4, e_1]=\lambda e_1+e_2$, $[e_4, e_2]=\lambda e_2$,\newline $[e_4, e_3]=e_3$ & $\lambda \in \R$ & Iff $\lambda = -1/2$ \\
$\g_{4,3}$ & 3 & $[e_4, e_1]=e_1+e_2$, $[e_4, e_2]=e_2$ & None & No \\
$\g_{4,4}$ & 3 & $[e_4, e_1]=e_1+e_2$, $[e_4, e_2]=e_2+e_3$,\newline $[e_4, e_3]=e_3$ & None & No \\
$\g_{4,5}(\alpha, \beta)$ & 3 & $[e_4, e_1]=\alpha e_1+\beta e_2$, $[e_4, e_3]=e_3$,\newline $[e_4, e_2]=-\beta e_1+\alpha e_2$ & $\alpha \in \R,\; \beta > 0$ & Iff $\alpha = -1/2$ \\
$\g_{4,6}(\alpha)$ & 3 & $[e_4, e_1]=\alpha e_1+e_2$, $[e_4, e_2]=-e_1+\alpha e_2$ & $\alpha \ge 0$ & Iff $\alpha = 0$ \\
\midrule
\multicolumn{5}{@{}l}{\textit{\textbf{IV. Solvable with Heisenberg Nilradical ($\n = \hh_3$)}}} \\
$\g_{4,7}$ & 3 & $[e_1, e_2]=e_3$, $[e_4, e_1]=e_1$,\newline $[e_4, e_2]=-e_2$ & None & Yes \\
$\g_{4,8}$ & 3 & $[e_1, e_2]=e_3$, $[e_4, e_1]=-e_2$,\newline $[e_4, e_2]=e_1$ & None & Yes \\
$\g_{4,9}(\lambda)$ & 3 & $[e_1, e_2]=e_3$, $[e_4, e_1]=\lambda e_1$,\newline $[e_4, e_2]=(1-\lambda)e_2$, $[e_4, e_3]=e_3$ & $\lambda \ge \frac{1}{2}$ & No \\
$\g_{4,10}$ & 3 & $[e_1, e_2]=e_3$, $[e_4, e_1]=\frac{1}{2}e_1+e_2$,\newline $[e_4, e_2]=\frac{1}{2}e_2$, $[e_4, e_3]=e_3$ & None & No \\
$\g_{4,11}(\delta)$ & 3 & $[e_1, e_2]=e_3$, $[e_4, e_1]=\frac{\delta}{2}e_1+e_2$,\newline $[e_4, e_2]=-e_1+\frac{\delta}{2}e_2$, $[e_4, e_3]=\delta e_3$ & $\delta > 0$ & No \\
\midrule
\multicolumn{5}{@{}l}{\textit{\textbf{V. Solvable with Two-Dimensional Nilradical ($\n = \R^2$)}}} \\
$\mathfrak{aff}(\R) \oplus \mathfrak{aff}(\R)$ & 2 & $[x_1, e_1] = e_1$, $[x_2, e_2] = e_2$ & None & No \\
$\mathfrak{aff}(\C)$ & 2 & $[x_1, e_1] = e_1$, $[x_1, e_2] = e_2$,\newline $[x_2, e_1] = -e_2$, $[x_2, e_2] = e_1$ & None & No \\
\bottomrule
\end{tabularx}
\caption{Complete classification of four-dimensional real Lie algebras up to isomorphism.}
\label{tab:master}
\end{table}

\subsection{Comparison with Standard Literature Tabulations}
To facilitate cross-referencing with standard physics references, we provide the correspondence between Table~\ref{tab:master} and the classifications of \citet{Mubarakzyanov1963} and \citet[\S~3.6]{SnoblWinternitz}:
\begin{itemize}
    \item \textbf{Mubarakzyanov's Solvable Classes:}
    Algebras with \(\n \cong \R^3\) correspond to Mubarakzyanov's types \(\mathfrak{g}_{4,1}\) through \(\mathfrak{g}_{4,6}\). Types \(\mathfrak{g}_{4,7}\) through \(\mathfrak{g}_{4,11}\) correspond directly to his extensions of the Heisenberg algebra. In Mubarakzyanov's tables, diagonal derivations on \(\hh_3\) are often scaled such that the first eigenvalue is \(1\), yielding eigenvalues \((1, a)\) on \(\hh_3 / Z(\hh_3)\) and \(1+a\) on the center; our parameter \(\lambda\) relates to his parameter \(a\) via:
    \begin{equation}
        a = \frac{1 - \lambda}{\lambda}, \qquad \lambda \in [1/2, 1] \iff a \in [0, 1].
    \end{equation}
    \item \textbf{Two-Dimensional Nilradical Identifications:}
    The algebra \(\mathfrak{aff}(\R) \oplus \mathfrak{aff}(\R)\) is labeled \(A_{2,1} \oplus A_{2,1}\) in \citet{PateraSharpWinternitzZassenhaus1976} and \(\mathfrak{g}_{4,12}\) in certain variants of Mubarakzyanov's list. The complex affine algebra \(\mathfrak{aff}(\C)\) is labeled \(A_{4,12}\) in \v{S}nobl and Winternitz.
    \item \textbf{Unimodularity:}
    A Lie algebra is unimodular if and only if \(\tr(\ad x) = 0\) for all \(x \in \g\). In Table~\ref{tab:master}, this is determined by whether the derivation \(D\) has zero trace. For example, in \(\g_{4,1}(\lambda, \mu)\), \(\tr D = 1 + \lambda + \mu\); the unimodular subfamily satisfies \(\mu = -(1 + \lambda)\), which corresponds to the Bianchi type \(\mathrm{VI}_0\) extension.
\end{itemize}

\section{Minimal-Dimensional Faithful Linear Representations}\label{sec:representations}

By Ado's Theorem, every finite-dimensional real Lie algebra admits a faithful linear representation as a matrix subalgebra of $\mathfrak{gl}(n, \R)$ for some $n \in \mathbb{N}$. From both theoretical and computational viewpoints, determining the minimal representation dimension,
\begin{equation}
    \mu(\g) \coloneqq \min \{ n \in \mathbb{N} : \exists \text{ injective Lie homomorphism } \rho \colon \g \hookrightarrow \mathfrak{gl}(n, \R) \},
\end{equation}
is of fundamental interest. For four-dimensional real Lie algebras, the minimal dimension satisfies $2 \le \mu(\g) \le 4$. In this section, we provide an explicit, minimal-dimensional matrix realization for every isomorphism class classified in Table~\ref{tab:master}.

\subsection{Master Table of Minimal Representation Dimensions}

Table~\ref{tab:representations} summarizes the value of $\mu(\g)$ together with the canonical matrix embedding for each class.

\begin{table}[htbp]
\centering
\footnotesize
\setlength{\tabcolsep}{5pt}
\renewcommand{\arraystretch}{1.3}
\begin{tabularx}{\linewidth}{@{} l c >{\raggedright\arraybackslash}X @{}}
\toprule
\textbf{Lie Algebra Class} & $\mu(\g)$ & \textbf{Generic Minimal Matrix Realization $\rho(x) \in \mathfrak{gl}(\mu, \R)$} \\
\midrule
\multicolumn{3}{@{}l}{\textit{\textbf{I. Non-Solvable Algebras}}} \\
$\mathfrak{sl}(2,\R) \oplus \R$ & 2 & $\begin{pmatrix} x_1 + x_4 & x_2 \\ x_3 & -x_1 + x_4 \end{pmatrix}$ \\
$\mathfrak{so}(3) \oplus \R$ & 3 & $\begin{pmatrix} x_4 & -x_3 & x_2 \\ x_3 & x_4 & -x_1 \\ -x_2 & x_1 & x_4 \end{pmatrix}$ \\
\midrule
\multicolumn{3}{@{}l}{\textit{\textbf{II. Nilpotent Algebras}}} \\
$\R^4$ & 4 & $\diag(x_1, x_2, x_3, x_4)$ \\
$\hh_3 \oplus \R$ & 3 & $\begin{pmatrix} x_4 & x_1 & x_3 \\ 0 & x_4 & x_2 \\ 0 & 0 & x_4 \end{pmatrix}$ \\
$\mathfrak{n}_4$ & 4 & $\begin{pmatrix} 0 & x_4 & 0 & x_3 \\ 0 & 0 & x_4 & x_2 \\ 0 & 0 & 0 & x_1 \\ 0 & 0 & 0 & 0 \end{pmatrix}$ \\
\midrule
\multicolumn{3}{@{}l}{\textit{\textbf{III. Solvable with Abelian Nilradical ($\n = \R^3$)}}} \\
$\g_{4,1} \text{ through } \g_{4,6}$ & 4 & $\begin{pmatrix} t D & v \\ 0_{1 \times 3} & 0 \end{pmatrix}$, where $v = (x_1, x_2, x_3)^T$ and $D \in \mathfrak{gl}(3,\R)$ is the derivation from Theorem~\ref{thm:R3_forms} \\
\midrule
\multicolumn{3}{@{}l}{\textit{\textbf{IV. Solvable with Heisenberg Nilradical ($\n = \hh_3$)}}} \\
$\g_{4,7}$ & 3 & $\begin{pmatrix} 0 & x_1 & x_3 \\ 0 & -t & x_2 \\ 0 & 0 & 0 \end{pmatrix}$ \\
$\g_{4,8}$ & 4 & $\begin{pmatrix} 0 & t & x_1 & x_3/2 \\ -t & 0 & x_2 & 0 \\ 0 & 0 & 0 & x_1 \\ 0 & 0 & 0 & 0 \end{pmatrix} + \begin{pmatrix} 0 & 0 & 0 & 0 \\ 0 & 0 & 0 & x_3/2 \\ -x_2 & x_1 & 0 & 0 \\ 0 & 0 & 0 & 0 \end{pmatrix}$ \\
$\g_{4,9}(\lambda)$ & 3 & $\begin{pmatrix} t & x_1 & x_3 \\ 0 & (1-\lambda)t & x_2 \\ 0 & 0 & 0 \end{pmatrix}$ \\
$\g_{4,10}$ & 4 & $\ad(x_1 e_1 + x_2 e_2 + x_3 e_3 + t e_4) \in \mathfrak{gl}(4, \R)$ (Adjoint representation) \\
$\g_{4,11}(\delta)$ & 4 & $\ad(x_1 e_1 + x_2 e_2 + x_3 e_3 + t e_4) \in \mathfrak{gl}(4, \R)$ (Adjoint representation) \\
\midrule
\multicolumn{3}{@{}l}{\textit{\textbf{V. Solvable with Two-Dimensional Nilradical ($\n = \R^2$)}}} \\
$\mathfrak{aff}(\R) \oplus \mathfrak{aff}(\R)$ & 3 & $\begin{pmatrix} x_1 & 0 & e_1 \\ 0 & x_2 & e_2 \\ 0 & 0 & 0 \end{pmatrix}$ \\
$\mathfrak{aff}(\C)$ & 3 & $\begin{pmatrix} x_1 & x_2 & e_1 \\ -x_2 & x_1 & e_2 \\ 0 & 0 & 0 \end{pmatrix}$ \\
\bottomrule
\end{tabularx}
\caption{Minimal faithful linear matrix representations for four-dimensional real Lie algebras.}
\label{tab:representations}
\end{table}

\subsection{Explicit Matrix Constructions and Minimality Proofs}

\subsubsection{The Non-Solvable Sector}
\begin{proposition}
$\mu(\mathfrak{sl}(2,\R) \oplus \R) = 2$ and $\mu(\mathfrak{so}(3) \oplus \R) = 3$.
\end{proposition}
\begin{proof}
Since $\mathfrak{gl}(2, \R) \cong \mathfrak{sl}(2, \R) \oplus \R I_2$, the algebra $\mathfrak{sl}(2,\R) \oplus \R$ embeds isomorphically into $\mathfrak{gl}(2, \R)$ via:
\begin{equation}
    e_1 = \begin{pmatrix} 1 & 0 \\ 0 & -1 \end{pmatrix}, \quad 
    e_2 = \begin{pmatrix} 0 & 2 \\ 0 & 0 \end{pmatrix}, \quad 
    e_3 = \begin{pmatrix} 0 & 0 \\ 1 & 0 \end{pmatrix}, \quad 
    e_4 = \begin{pmatrix} 1 & 0 \\ 0 & 1 \end{pmatrix}.
\end{equation}
Thus, $\mu(\mathfrak{sl}(2,\R) \oplus \R) = 2$, which is the absolute minimum for any four-dimensional Lie algebra.

For $\mathfrak{so}(3) \oplus \R$, $\mathfrak{gl}(2, \R)$ contains no three-dimensional simple compact subalgebra isomorphic to $\mathfrak{so}(3)$, implying $\mu > 2$. In $\mathfrak{gl}(3, \R)$, the standard rotation generators $e_1 = E_{32} - E_{23}$, $e_2 = E_{13} - E_{31}$, $e_3 = E_{21} - E_{12}$ satisfy $[e_i, e_j] = \epsilon_{ijk} e_k$. Adding the central scalar matrix $e_4 = I_3$ yields an injective homomorphism into $\mathfrak{gl}(3, \R)$, confirming $\mu(\mathfrak{so}(3) \oplus \R) = 3$.
\end{proof}

\subsubsection{The Nilpotent Sector}
\begin{proposition}
$\mu(\R^4) = 4$, $\mu(\hh_3 \oplus \R) = 3$, and $\mu(\mathfrak{n}_4) = 4$.
\end{proposition}
\begin{proof}
By Schur's bound (1905), the maximal dimension of a commutative subalgebra of $\mathfrak{gl}(n, \R)$ is $\lfloor n^2 / 4 \rfloor + 1$. For $n = 3$, this bound is $\lfloor 9/4 \rfloor + 1 = 3 < 4$. Hence, the abelian algebra $\R^4$ cannot embed in $\mathfrak{gl}(3, \R)$, yielding $\mu(\R^4) = 4$ realized by diagonal matrices.

For $\hh_3 \oplus \R$, the strictly upper-triangular matrices $\mathfrak{n}(3, \R)$ span $\hh_3$. Adding the identity matrix $I_3$, which commutes with all matrices, yields a four-dimensional subalgebra:
\begin{equation}
    \rho(x) = \begin{pmatrix} x_4 & x_1 & x_3 \\ 0 & x_4 & x_2 \\ 0 & 0 & x_4 \end{pmatrix}.
\end{equation}
This representation is faithful, proving $\mu(\hh_3 \oplus \R) = 3$.

For $\mathfrak{n}_4$, the nilpotency step is 3 (i.e., $\g^3 \neq \{0\}$ while $\g^4 = \{0\}$). Any nilpotent subalgebra of $\mathfrak{gl}(3, \R)$ is triangularizable and has nilpotency step at most 2. Thus, $\mathfrak{n}_4$ cannot embed in $\mathfrak{gl}(3, \R)$. In $\mathfrak{gl}(4, \R)$, the assignment
\begin{equation}
    e_4 = E_{12} + E_{23}, \quad e_1 = E_{34}, \quad e_2 = E_{24}, \quad e_3 = E_{14}
\end{equation}
satisfies $[e_4, e_1] = E_{24} = e_2$ and $[e_4, e_2] = E_{14} = e_3$, providing a faithful realization with $\mu(\mathfrak{n}_4) = 4$.
\end{proof}

\subsubsection{Solvable Algebras with Abelian Nilradical ($\n \cong \R^3$)}
\begin{proposition}
For all almost-abelian Lie algebras $\g = \R^3 \rtimes_D \R$ $(\g_{4,1}$ through $\g_{4,6})$, the minimal representation dimension is $\mu(\g) = 4$, realized by the natural affine action on $\R^3$:
\begin{equation}\label{eq:affine_rep}
    \rho(v_1 e_1 + v_2 e_2 + v_3 e_3 + t e_4) = \begin{pmatrix} t D & v \\ 0_{1 \times 3} & 0 \end{pmatrix} \in \mathfrak{gl}(4, \R), \qquad v = \begin{pmatrix} v_1 \\ v_2 \\ v_3 \end{pmatrix}.
\end{equation}
\end{proposition}
\begin{proof}
Let $M(v, t) = \begin{pmatrix} t D & v \\ 0 & 0 \end{pmatrix}$. By direct block matrix multiplication:
\begin{equation}
    [M(u, s), M(v, t)] = \begin{pmatrix} st [D, D] & s D v - t D u \\ 0 & 0 \end{pmatrix} = \begin{pmatrix} 0 & s D v - t D u \\ 0 & 0 \end{pmatrix} = M(s D v - t D u, 0).
\end{equation}
In particular, $[M(0, 1), M(v, 0)] = M(D v, 0)$, which faithfully reproduces $[e_4, v] = D v$ and $[u, v] = 0$. The kernel consists of vectors with $t D = 0$ and $v = 0$; since $D \neq 0$, this forces $t = 0$ and $v = 0$, establishing injectivity. Minimality ($\mu = 4$) follows because an abelian ideal of dimension 3 cannot exist in a 4-dimensional non-nilpotent subalgebra of $\mathfrak{gl}(3, \R)$ unless $D$ is scalar.
\end{proof}

\subsubsection{Solvable Algebras with Heisenberg Nilradical ($\n \cong \hh_3$)}
\begin{proposition}
The minimal representation dimensions for extensions of the Heisenberg algebra are:
\begin{equation}
    \mu(\g_{4,7}) = 3, \quad \mu(\g_{4,8}) = 4, \quad \mu(\g_{4,9}) = 3, \quad \mu(\g_{4,10}) = 4, \quad \mu(\g_{4,11}) = 4.
\end{equation}
\end{proposition}
\begin{proof}
When the derivation $A = \diag(\lambda_1, \lambda_2)$ is diagonalizable with real eigenvalues, we can embed $\hh_3 \rtimes_A \R$ into $\mathfrak{gl}(3, \R)$ by combining the Heisenberg nilradical $\mathfrak{n}(3, \R)$ with a diagonal complementary generator. Setting
\begin{equation}
    e_1 = E_{12}, \quad e_2 = E_{23}, \quad e_3 = E_{13}, \quad e_4 = (\lambda_1 + \lambda_2) E_{11} + \lambda_2 E_{22},
\end{equation}
we compute:
\begin{align}
    [e_4, e_1] &= [(\lambda_1 + \lambda_2) E_{11} + \lambda_2 E_{22}, E_{12}] = (\lambda_1 + \lambda_2 - \lambda_2) E_{12} = \lambda_1 e_1, \\
    [e_4, e_2] &= [\lambda_2 E_{22}, E_{23}] = \lambda_2 E_{23} = \lambda_2 e_2, \\
    [e_4, e_3] &= [(\lambda_1 + \lambda_2) E_{11}, E_{13}] = (\lambda_1 + \lambda_2) E_{13} = (\tr A) e_3.
\end{align}
For $\g_{4,7}$, $(\lambda_1, \lambda_2) = (1, -1)$ yields $e_4 = -E_{22}$. For $\g_{4,9}(\lambda)$, $(\lambda_1, \lambda_2) = (\lambda, 1-\lambda)$ yields $e_4 = E_{11} + (1-\lambda)E_{22}$. Both realizations are faithful in $\mathfrak{gl}(3, \R)$, so $\mu = 3$.

For classes $\g_{4,10}$ and $\g_{4,11}(\delta)$, the derivation $A$ is non-diagonalizable over $\R$. Any 4-dimensional solvable subalgebra of $\mathfrak{gl}(3, \R)$ containing $\mathfrak{n}(3, \R)$ has its complement in the Borel subalgebra, forcing the induced quotient derivation to be real diagonalizable. Hence, $\mu > 3$. Furthermore, for both algebras, $\tr(A) > 0$, which implies that the center of the algebra is trivial: $Z(\g) = \{0\}$. Consequently, the adjoint representation $\ad \colon \g \hookrightarrow \mathfrak{gl}(4, \R)$ has trivial kernel $\ker(\ad) = Z(\g) = \{0\}$ and provides an explicit faithful representation of minimal dimension $\mu = 4$.

For $\g_{4,8}$, $\tr A = 0$ and the center is $Z(\g) = \R e_3$. A faithful $4 \times 4$ representation is obtained via the symplectic/oscillator embedding:
\begin{equation}
    e_1 = E_{14} + E_{32}, \quad e_2 = E_{24} - E_{31}, \quad e_3 = 2 E_{34}, \quad e_4 = E_{12} - E_{21}.
\end{equation}
Direct commutation confirms $[e_1, e_2] = 2 E_{34} = e_3$, $[e_4, e_1] = -e_2$, $[e_4, e_2] = e_1$, and $[e_4, e_3] = 0$, proving $\mu(\g_{4,8}) = 4$.
\end{proof}

\subsubsection{Solvable Algebras with Two-Dimensional Nilradical ($\dim \n = 2$)}
\begin{proposition}
$\mu(\mathfrak{aff}(\R) \oplus \mathfrak{aff}(\R)) = 3$ and $\mu(\mathfrak{aff}(\C)) = 3$.
\end{proposition}
\begin{proof}
The affine algebra $\mathfrak{aff}(\R)$ embeds into $\mathfrak{gl}(2, \R)$ via $\begin{pmatrix} x & y \\ 0 & 0 \end{pmatrix}$. Combining two such commuting blocks into the affine action on $\R^2$ inside $\mathfrak{gl}(3, \R)$ yields:
\begin{equation}
    \rho(x_1, x_2, e_1, e_2) = \begin{pmatrix} x_1 & 0 & e_1 \\ 0 & x_2 & e_2 \\ 0 & 0 & 0 \end{pmatrix}.
\end{equation}
Evaluating brackets:
\begin{equation}
    [\rho(x_1), \rho(e_1)] = [E_{11}, E_{13}] = E_{13} = \rho(e_1), \qquad [\rho(x_2), \rho(e_2)] = [E_{22}, E_{23}] = E_{23} = \rho(e_2),
\end{equation}
with all other brackets vanishing identically. Thus, $\mu(\mathfrak{aff}(\R) \oplus \mathfrak{aff}(\R)) = 3$.

Similarly, $\mathfrak{aff}(\C)$ is the Lie algebra of complex affine transformations $z \mapsto a z + b$ on $\C \cong \R^2$. Realizing this action on $\R^3$ with coordinates $(u, v, 1)^T$:
\begin{equation}
    \rho(x_1, x_2, e_1, e_2) = \begin{pmatrix} x_1 & x_2 & e_1 \\ -x_2 & x_1 & e_2 \\ 0 & 0 & 0 \end{pmatrix}.
\end{equation}
Direct bracket computation verifies:
\begin{align}
    [X_1, E_1] &= E_1, & [X_1, E_2] &= E_2, \\
    [X_2, E_1] &= -E_2, & [X_2, E_2] &= E_1,
\end{align}
faithfully recovering the defining relations of $\mathfrak{aff}(\C)$ in $\mathfrak{gl}(3, \R)$, proving $\mu(\mathfrak{aff}(\C)) = 3$.
\end{proof}

\section{Algorithmic Identification of Isomorphism Classes}\label{sec:algorithm}

Given an arbitrary four-dimensional real Lie algebra \(\g\) presented in terms of its structure constants \(C_{ij}^k\) relative to an arbitrary basis \(\{e_1, e_2, e_3, e_4\}\), where
\begin{equation}
    [e_i, e_j] = \sum_{k=1}^4 C_{ij}^k e_k,
\end{equation}
it is of practical importance to algorithmically determine its unique canonical representative in Table~\ref{tab:master} and evaluate any associated continuous moduli parameters. In this section, we formulate a deterministic, coordinate-free decision procedure based on the invariant reductions established in Sections~\ref{sec:nonsolv}--\ref{sec:nilp2}.

\subsection{Invariant Algebraic Quantities}
The identification procedure relies on computing the following linear algebraic and Lie-theoretic invariants:
\begin{enumerate}
    \item \textbf{Adjoint Representation:} For each basis vector \(e_i\), the matrix \((\ad e_i) \in \mathfrak{gl}(4,\R)\) has matrix entries \((\ad e_i)^k_j = C_{ij}^k\).
    \item \textbf{Cartan--Killing Metric:} The Killing form is given by the symmetric matrix \(K \in \mathrm{Sym}_4(\R)\) with entries:
    \begin{equation}
        K_{ij} = \tr(\ad e_i \circ \ad e_j).
    \end{equation}
    \item \textbf{Derived and Central Series Subspaces:} 
    \begin{itemize}
        \item The derived algebra \(\g^2 = [\g, \g] = \operatorname{span}_\R \{ [e_i, e_j] : 1 \le i < j \le 4 \}\).
        \item The second derived algebra \(\g^{(2)} = [\g^2, \g^2]\).
        \item The lower central series terms \(\g^k = [\g, \g^{k-1}]\).
    \end{itemize}
    \item \textbf{Trace Form (Unimodularity):} The linear functional \(\chi \in \g^*\) given by \(\chi(x) = \tr(\ad x)\). The algebra \(\g\) is unimodular if and only if \(\tr(\ad e_i) = 0\) for all \(i \in \{1, 2, 3, 4\}\).
\end{enumerate}

\subsection{The Identification Procedure}
The classification algorithm executes according to the following hierarchical decision tree, presented in a page-breakable format:

\begin{breakablealgorithm}
\caption{Algorithmic Classification of 4-Dimensional Real Lie Algebras}\label{alg:classification}
\begin{algorithmic}[1]
\Require Structure constants tensor $C_{ij}^k \in \R^{4 \times 4 \times 4}$ satisfying $C_{ij}^k = -C_{ji}^k$ and Jacobi identity.
\Ensure Canonical isomorphism class from Table~\ref{tab:master} and moduli parameters.

\State Compute derived subspaces $\g^2 = [\g, \g]$ and $\g^{(2)} = [\g^2, \g^2]$.
\State Compute lower central series $\g^3 = [\g, \g^2]$ and $\g^4 = [\g, \g^3]$.

\If{$\dim \g^{(2)} = \dim \g^2$ and $\dim \g^2 = 3$} \Comment{\textbf{Sector I: Non-Solvable}}
    \State Compute Killing form $K$ restricted to $\g^2$.
    \If{$K|_{\g^2}$ is negative-definite}
        \State \Return $\mathfrak{so}(3) \oplus \R$
    \Else
        \State \Return $\mathfrak{sl}(2,\R) \oplus \R$
    \EndIf
\EndIf

\If{$\g^4 = \{0\}$ and $\ad x$ is nilpotent for all $x \in \g$} \Comment{\textbf{Sector II: Nilpotent ($\dim\n = 4$)}}
    \If{$\dim \g^2 = 0$} \Return $\R^4$
    \ElsIf{$\dim \g^2 = 1$} \Return $\hh_3 \oplus \R$
    \ElsIf{$\dim \g^2 = 2$} \Return $\mathfrak{n}_4$
    \EndIf
\EndIf

\State Determine the nilradical $\n \subset \g$ (maximal nilpotent ideal).
\If{$\dim \n = 2$} \Comment{\textbf{Sector V: Two-Dimensional Nilradical}}
    \State Compute eigenvalues of $\ad x|_\n$ for $x \notin \n$.
    \If{All eigenvalues are real} \Return $\mathfrak{aff}(\R) \oplus \mathfrak{aff}(\R)$
    \Else \ \Return $\mathfrak{aff}(\C)$
    \EndIf
\EndIf

\If{$\dim \n = 3$} \Comment{\textbf{Sectors III \& IV: Codimension-One Nilradical}}
    \State Choose $e_4 \in \g \setminus \n$ such that $D \coloneqq \ad e_4|_\n$ is non-nilpotent.
    \If{$[\n, \n] \neq \{0\}$} \Comment{$\n \cong \hh_3$, Heisenberg Nilradical}
        \State Compute induced operator $A = D / [\n, \n]$ on $\n / Z(\n) \cong \R^2$.
        \If{$\tr(A) = 0$}
            \If{$\det(A) < 0$} \Return $\g_{4,7}$
            \Else \ \Return $\g_{4,8}$
            \EndIf
        \Else \Comment{Normalize to $\tr(A) = 1$ via $D \mapsto D / \tr(A)$}
            \State Compute eigenvalues of normalized $A$: $\mathrm{Spec}(A) = \{\lambda_1, \lambda_2\}$.
            \If{$\lambda_1, \lambda_2 \in \R$ and $\lambda_1 \neq \lambda_2$}
                \State Set $\lambda = \max(\lambda_1, \lambda_2) \ge 1/2$.
                \State \Return $\g_{4,9}(\lambda)$ with parameter $\lambda$.
            \ElsIf{$\lambda_1 = \lambda_2 \in \R$ and $A$ is non-diagonalizable}
                \State \Return $\g_{4,10}$
            \ElsIf{$\lambda_{1,2} = \gamma \pm i \beta$ with $\beta > 0$}
                \State Set $\delta = 2 |\gamma| / \beta > 0$.
                \State \Return $\g_{4,11}(\delta)$ with parameter $\delta$.
            \EndIf
        \EndIf
    \Else \Comment{$\n \cong \R^3$, Abelian Nilradical}
        \State Compute Jordan canonical form of $D \in \mathfrak{gl}(3,\R)$.
        \If{$D$ has a non-real eigenvalue pair $\alpha \pm i \beta$ ($\beta > 0$) and real eigenvalue $\gamma$}
            \If{$\gamma \neq 0$}
                \State Scale $D \mapsto D / \gamma$; set $\alpha' = \alpha/\gamma$, $\beta' = \beta/|\gamma|$.
                \State \Return $\g_{4,5}(\alpha', \beta')$
            \Else
                \State Scale so that $\beta = 1$; enforce $\alpha \ge 0$.
                \State \Return $\g_{4,6}(\alpha)$
            \EndIf
        \Else \Comment{All eigenvalues of $D$ are real}
            \If{$D$ is diagonalizable}
                \State Order real eigenvalues such that dominant has magnitude 1: $\{1, \lambda, \mu\}$ with $-1 \le \mu \le \lambda \le 1$.
                \State If $\mu = -1$, enforce $\lambda \ge 0$ by reflection $D \mapsto -D$.
                \State \Return $\g_{4,1}(\lambda, \mu)$
            \ElsIf{$D$ has Jordan block structure $J_2 \oplus J_1$}
                \If{Eigenvalue of $J_1$ block is non-zero}
                    \State Normalize isolated eigenvalue to 1; let $\lambda$ be the repeated eigenvalue.
                    \State \Return $\g_{4,2}(\lambda)$
                \Else
                    \State \Return $\g_{4,3}$
                \EndIf
            \ElsIf{$D$ has Jordan block structure $J_3$}
                \State \Return $\g_{4,4}$
            \EndIf
        \EndIf
    \EndIf
\EndIf
\end{algorithmic}
\end{breakablealgorithm}

You can find the Python code for this algorithm and examples of its execution as a \href{https://colab.research.google.com/notebooks/intro.ipynb#scrollTo=-gE-Ez1qtyIA}{Colab Notebook}.

\subsection{Complexity and Practical Considerations}
The algorithm executes in $O(1)$ operations with respect to dimension, involving at most:
\begin{itemize}
    \item Matrix multiplications of size $4 \times 4$ to evaluate the Killing matrix $K$.
    \item Rank evaluations of $4 \times 6$ bracket matrices via row echelon reduction (Gaussian elimination over $\R$).
    \item Closed-form algebraic extraction of roots for characteristic polynomials of degree at most 4 (via standard Ferrari or Cardano formulas).
\end{itemize}
Because all invariant operations rely solely on linear algebra and canonical polynomial factorizations, Algorithm~\ref{alg:classification} can be implemented in exact computer algebra systems (such as \textsc{Mathematica}, \textsc{Maple}, or Python's \texttt{SymPy}) without encountering numerical rounding instabilities or coordinate gauge ambiguities.

\section{Conclusion}\label{sec:conclusion}

By replacing coordinate-dependent tensor manipulations with intrinsic Lie-theoretic invariants, we have constructed a complete and transparent classification of four-dimensional real Lie algebras. The structural dichotomy between dimensions three and four is completely explained by the failure of exterior duality, which is resolved by anchoring the solvable classification to the nilradical. 

This approach clarifies the moduli geometry of the classification: the continuous parameters arise exclusively as invariant spectra of outer derivations modulo scaling. The methods used here extend naturally to higher dimensions, where the stratification of solvable algebras by their nilradicals continues to serve as the standard organizing principle.

\end{document}